\newcommand{\Dim}{\mathrm{dim}}
\newcommand{\gR}{g_{\mathrm R}}

\newcommand{\R}{\mathds R}
\newcommand{\C}{\mathds C}
\newcommand{\N}{\mathds N}

\newcommand{\Iso}{\mathrm{Iso}}

\documentclass[twoside,final,a4paper,10pt]{amsart}

\usepackage{times}
\usepackage{dsfont}

\numberwithin{equation}{section}

\title[Periodic geodesics in compact stationary manifolds]{Periodic geodesics and geometry of compact Lorentzian
manifolds with a Killing vector field}
\author[J. L. Flores]{Jos\'{e} Luis Flores}
\address{Departamento de \'{A}lgebra, Geometr\'{\i}a y Topolog\'{\i}a.\hfill\break\indent
 Facultad de Ciencias, Universidad de M\'{a}laga.\hfill\break\indent
 Campus Teatinos s/n, 29071 M\'{a}laga, Spain}
 \email{floresj@agt.cie.uma.es}
\author[M. \'{A}ngel Javaloyes]{Miguel \'{A}ngel Javaloyes}
\address{Departamento de Geometr\'{\i}a y Topolog\'{\i}a.\hfill\break\indent
 Facultad de Ciencias, Universidad de Granada.\hfill\break\indent
 Campus Fuentenueva s/n, 18071 Granada, Spain}
 \email{ma.javaloyes@gmail.com}
\author[P.\ Piccione]{Paolo Piccione}
\address{Departamento de Matem\'atica,\hfill\break\indent
Universidade de S\~ao Paulo, \hfill\break\indent Rua do Mat\~ao
1010,\hfill\break\indent CEP 05508-900, S\~ao Paulo, SP, Brazil}
\email{piccione.p@gmail.com}
\curraddr{Department of Mathematics, \hfill\break\indent
University of Murcia, Campus de Espinardo\hfill\break\indent
30100 Espinardo, Murcia, \hfill\break\indent Spain}
\thanks{The authors are grateful to Andrea Spiro for giving several interesting suggestions
on isometric group actions on manifolds. The first author thanks
the Departamento de Matem\'atica, Universidade de S\~ao Paulo,
where this work started, for its kind hospitality during his stay
there. He is partially supported by Spanish MEC-FEDER Grant
MTM2007-60731 and Regional J. Andaluc\'{i}a Grant P06-FQM-01951.
The second author was partially supported by Regional J.
Andaluc\'{\i}a Grant P06-FQM-01951 and by Spanish MEC Grant
MTM2007-64504. The third author is sponsored by Capes, Brasil,
Grant BEX 1509-08-0.}


\date{February 25th, 2009}

\begin{document}


\theoremstyle{plain}\newtheorem*{teon}{Theorem}
\theoremstyle{definition}\newtheorem*{defin*}{Definition}
\theoremstyle{plain}\newtheorem{teo}{Theorem}[section]
\theoremstyle{plain}\newtheorem{prop}[teo]{Proposition}
\theoremstyle{plain}\newtheorem{lem}[teo]{Lemma}
\theoremstyle{plain}\newtheorem{cor}[teo]{Corollary}
\theoremstyle{definition}\newtheorem{defin}[teo]{Definition}
\theoremstyle{remark}\newtheorem{rem}[teo]{Remark}
\theoremstyle{plain} \newtheorem{assum}[teo]{Assumption}
\theoremstyle{definition}\newtheorem{example}{Example}
\swapnumbers
\theoremstyle{plain}
\newtheorem*{acknowledgement}{Acknowledgements}
\theoremstyle{definition}\newtheorem*{notation}{Notation}


\begin{abstract} We study the geometry and the periodic geodesics of
a compact Lorentzian manifold that has a Killing vector field
which is timelike somewhere. Using a compactness argument for
subgroups of the isometry group, we prove the existence of one
timelike non self-intersecting periodic geodesic. If the Killing
vector field is never vanishing, then there are at least two
distinct periodic geodesics; as a special case, compact stationary
manifolds have at least two periodic timelike geodesics. We also
discuss some properties of the topology of such manifolds. In
particular, we show that a compact manifold $M$ admits a
Lorentzian metric with a never vanishing Killing  vector field which is
timelike somewhere if and only if $M$ admits a smooth circle
action without fixed points.
\end{abstract}

\maketitle
\begin{section}{Introduction}
A classical subject in Geometry is the question of existence and
multiplicity of periodic geodesics.  A well known result in
differential geometry establishes the existence of a periodic
geodesic in every compact Riemannian manifold; a Lorentzian analog
of this result is still an open problem in its full generality.
However, there are remarkable partial results in this direction.
An earlier result by Tipler (see \cite{Tipler}) gives the
existence of one periodic timelike geodesic in compact Lorentzian
manifolds that admit a regular covering which has a compact Cauchy
surface. Recently, this result has been extended by Guediri
\cite{Gue1, Gue2, Gue5} and S\'anchez \cite{San2} to the case that
the Cauchy surface in the covering is not necessarily compact, but
assuming certain hypotheses on the group of deck transformations.
The existence of a periodic timelike geodesic has been established
also by Galloway in \cite{Gal1}, where he proves the existence of
a longest periodic timelike curve, which is necessarily a
geodesic, in each \emph{stable} free timelike homotopy class.   In
\cite{Gal2}, the same author proves the existence of a causal
(i.e., nonspacelike) periodic geodesic in any compact
two-dimensional Lorentzian manifold. More recently, Guediri
 has proved that compact flat spacetimes contain a
causal periodic geodesic \cite{Gue1}, and that such spacetimes
contain a periodic timelike geodesic if and only if the
fundamental group of the underlying manifold contains a nontrivial
timelike translation \cite{Gue2}. Non existence results for
periodic causal geodesics are also available, see \cite{Gal2,
Gue3, Gue4}.

In the particular case of \emph{static} compact Lorentzian
manifolds, using variational methods it has been established the
existence of a periodic geodesic in each free homotopy class
corresponding to an element of the fundamental group having finite
conjugacy class, \cite{CapMasPic}. This result has been
generalized to \emph{any} free homotopy class containing a
periodic timelike curve in \cite{San2}, obtaining in particular
that any compact static spacetime admits a periodic timelike
geodesic. As to  the stationary case, there exist in literature
some previous results when the spacetime admits a standard
stationary expression and, as a consequence, it is never compact
(see \cite{BilMerPic, Masiello2}). As suggested by the authors in
\cite{CapMasPic,Gue5,San2}, an interesting open question would be
to determine if existence results for periodic geodesics also hold
for any {\em stationary} compact Lorentzian manifold. Moreover, in
general there does not exist a globally hyperbolic covering for
this class of manifolds (see \cite[pag.~23]{San2}, observe also
that, unlike the static case, compact stationary Lorentzian
manifolds may be simply connected) and a different approach from
that of \cite{CapMasPic,Gue1,Gue2,Gue5,San2,Tipler} is needed.

In this paper, we answer positively to the above question by using the strong
relation between  Killing vector fields and geodesics, and a compactness criterion
for subgroups of the isometry group of a Lorentzian manifold. More precisely, we prove the following:
\begin{teon}
Let $(M,g)$ be a compact Lorentzian manifold with $\Dim(M)\ge2$
that admits a Killing vector field $K$ that is timelike somewhere.
Then there is some non trivial periodic non self-intersecting
timelike geodesic in $(M,g)$. If either one of the following two conditions
is satisfied, then there are at least two non trivial periodic non
self-intersecting geodesics in $M$:
\begin{itemize}
\item[(a)] $\max\limits_{q\in M}g(K_q,K_q)\ne0$;
\smallskip

\item[(b)] $K$ is never vanishing.
\end{itemize}
When either condition is satisfied, if in addition $K$ has at most one periodic integral line,
then there are infinitely many geometrically distinct non trivial periodic non self-intersecting
geodesics in $(M,g)$.
\end{teon}
Recall that two periodic geodesics $c_1,c_2:\R\to M$ are
geometrically distinct if the sets $c_1(\R)$ and $c_2(\R)$ are
distinct, i.e., if one cannot be obtained as an iteration of the
other.  The proof of our theorem employs Lie group techniques; the key result
is a compactness criterion for subgroups of the isometry group of a compact
Lorentzian manifold (Proposition~\ref{thm:elem}).
Using this criterion, one shows that a compact Lorentzian manifold with a Killing vector field
that is timelike somewhere also has a Killing vector field $K$
that is timelike somewhere, and all of whose integral curves are
periodic (this is called a \emph{closed} Killing  vector field in the
paper). In this situation, the periodic geodesics are given by
those integral curves of $K$ that pass through critical points for
the function $f=g(K,K)$ on $M$ (Lemma~\ref{thm:intcurvgeo}). The
existence of a never vanishing closed Killing vector field also
gives us information on the topology of this manifold, that can be
characterized as a compact manifold admitting a smooth action of
the circle $\mathds S^1$ all of whose orbits have finite
stabilizer or, equivalently, without fixed points
(Section~\ref{sec:topology}). Fixed point free actions of the
circle on low dimensional compact (simply connected) manifolds are
classified, see \cite{Sco, Sei} for actions in dimension $3$,
\cite{Fin1, Fin2} in dimension $4$ and \cite{Kol} in dimension
$5$.

The orbit space $M/\mathds S^1$ may fail to be a manifold when the
action is not free (see Remark~\ref{rem:kleinbottle}), but in
general it has the structure of a compact orbifold. Using
equivariant Ljusternik--Schnirelman category theory, one obtains a
slightly better estimate on the number of integral lines of $K$
that are geodesics (Corollary~\ref{thm:LS}); such number is
greater than or equal to the equivariant category
$\mathrm{cat}_{\mathds S^1}(M)$, which is greater than or equal to
$2$. If one assumes that
$(M,g)$ has a non closed Killing vector field that is timelike
somewhere, and that $\max g(K,K)$ is non zero,
then we prove that $(M,g)$ has infinitely many closed Killing
vector fields that are timelike somewhere, each of which gives
rise to (at least) two periodic geodesics. Under the further
assumption that the Killing vector field has at most one periodic
integral curve, then we prove that the family of periodic
geodesics produced in this way contains indeed infinitely many
geometrically distinct periodic geodesics. Explicit examples of
compact stationary Lorentzian manifolds with a timelike Killing
vector field all of whose integral lines \emph{but one} are non
periodic are given in Subsection~\ref{sub:someexamples}.

All the results of this paper apply in particular
to compact stationary manifolds, i.e., Lorentzian manifolds $(M,g)$ that admit an
\emph{everywhere timelike} Killing vector field $K$. It should be observed that,
in this case, there is an alternative approach to the problem which uses an auxiliary
Riemannian metric $\gR$ naturally associated to $g$ and $K$ (see \eqref{eq:LorRiem}).

A very interesting open question is whether a simply connected compact stationary
Lorentzian manifold has compact isometry group. The answer is yes in the real analytic
case (see \cite{Da}). In Subsection~\ref{sub:isometrygroup} we discuss briefly this issue,
giving some partial results towards a positive answer to the compactness question.
The question is studied in full generality in \cite{PicZeg}.

Finally, in Section~\ref{sec:final} we present a few related results on the existence
of periodic geodesics in compact semi-Riemannian manifolds, obtained using the techniques
discussed in the paper.
\end{section}

\begin{section}{Proof of the Theorem}\label{s2}
\label{sec:proofoftheorem}
\subsection{Preliminaries on Killing vector fields}

Semi-Riemannian manifolds are assumed to be connected with
dimension $n\geq 2$. In the case
of Lorentzian manifolds, we assume they are time-oriented (thus,
time-orientable). Our notation and conventions follow the standard
ones in Lorentzian Geometry, see \cite{BEE, ON, SW}.

Recall that a Killing vector field $K$ in a semi-Riemannian manifold
$(M,g)$ is a vector field whose flow preserves $g$, or,
equivalently, such that $\nabla K$ is skew--symmetric, i.e., the
bilinear map $T_pM\times T_pM\ni(v,w)\mapsto g(\nabla_vK,w)\in\R$
is skew--symmetric for all $p\in M$. Here $\nabla$ is the
Levi--Civita connection of $g$. If $K$ is a Killing vector field, then
the function $f=g(K,K)$ is constant along each integral curve of
$K$, namely, $K(f)=2g\big(\nabla_KK,K)=0$. The Lie bracket of
Killing vector fields on $M$ is a Killing vector field, so that
the space $\mathrm{Kill}(M,g)$ of all Killing vector fields on $M$
is a Lie algebra. Given a diffeomorphism $\Phi:M\to M$ and a
vector field $K$ on $M$, the \emph{push-forward} of $K$ by $\Phi$
is the vector field $\Phi_*(K)$ on $M$ defined by
$\Phi_*(K)_p=\mathrm d\Phi\big(\Phi^{-1}(p)\big)K_{\Phi^{-1}(p)}$
for all $p\in M$.

Let $(M,g)$ be a semi-Riemannian manifold, we will denote by
$\mathrm{Iso}(M,g)$ the Lie group of all isometries of $(M,g)$,
and by $\mathfrak{Iso}(M,g)$ its Lie algebra. It is well known
that the isometry group of a compact Riemannian manifold is a
compact Lie group. The natural action of $\mathrm{Iso}(M,g)$ on
$M$ is smooth, and the non trivial elements of $\mathrm{Iso}(M,g)$
cannot act trivially on a non empty open subset of $M$ (see \cite[Proposition 3.62]{ON}).

A correspondence between elements of $\mathfrak{Iso}(M,g)$ and
Killing vector fields on $M$ is obtained as follows. For $p\in M$
denote by $\beta_p:\mathrm{Iso}(M,g)\to M$ the smooth map
$\beta_p(\Phi)=\Phi(p)$. A Lie algebra anti-isomorphism
$\mathfrak{Iso}(M,g)\ni\mathfrak x\mapsto K^\mathfrak
x\in\mathrm{Kill}(M,g)$ is obtained by setting:
\begin{equation}\label{KillingLie}
K^\mathfrak x_p=\mathrm d\beta_p(1)\mathfrak x,\quad p\in M.
\end{equation}
For $\Phi\in\mathrm{Iso}(M,g)$ consider the Lie group isomorphism
$I_\Phi:\mathrm{Iso}(M,g)\to\mathrm{Iso}(M,g)$ given by
$I_\Phi(\Psi)=\Phi\circ\Psi\circ\Phi^{-1}$, and denote by
$\mathrm{Ad}_\Phi:\mathfrak{Iso}(M,g) \to\mathfrak{Iso}(M,g)$ its
differential at the identity. Using these notations, one has the
following immediate equality:
\begin{equation}\label{eq:compgammabeta}
\Phi\circ\beta_{\Phi^{-1}(p)}=\beta_p\circ I_\Phi,
\end{equation}
for all $\Phi\in\mathrm{Iso}(M,g)$ and all $p\in M$.
\begin{lem}\label{thm:Adpushforward}
For all $\Phi\in\mathrm{Iso}(M,g)$ and all $\mathfrak
x\in\mathfrak{Iso}(M,g)$, the following formula holds:
\[\Phi_*K^\mathfrak x=K^{\mathrm{Ad}_\Phi(\mathfrak x)}.\]
\end{lem}
\begin{proof}
A direct computation:
\begin{multline*}
\big(\Phi_*K^\mathfrak x\big)_p=\mathrm
d\Phi\big(\Phi^{-1}(p)\big)K^{\mathfrak x}_{\Phi^{-1}(p)}= \mathrm
d\Phi\big(\Phi^{-1}(p)\big)\,\mathrm
d\beta_{\Phi^{-1}(p)}(1)\mathfrak x\\= \mathrm
d\big(\Phi\circ\beta_{\Phi^{-1}(p)}\big)(1)\mathfrak x
\stackrel{\text{by~\eqref{eq:compgammabeta}}}=\mathrm
d\big(\beta_p\circ I_\Phi\big)(1)\mathfrak x= \mathrm
d\beta_p(1)\,\mathrm d I_\Phi(1)\mathfrak
x=K^{\mathrm{Ad}_\Phi(\mathfrak x)}_p.\qedhere
\end{multline*}

If $K$ is a Killing vector field for $(M,g)$ corresponding to the
element $\mathfrak x\in\mathfrak{Iso}(M,g)$, we will say that $K$
is \emph{closed} if the $1$-parameter group of isometries
$\big\{\exp(t\mathfrak x):t\in\R\big\}$ of $\mathrm{Iso}(M,g)$
generated by $\mathfrak x$ is closed. Here $\exp$ is the
exponential map of the Lie group $\mathrm{Iso}(M,g)$. The integral
curves of $K$ through some $p\in M$ of $K$ are given by
$t\mapsto\exp(t\mathfrak x)\cdot p$ (see \cite[Lemma 34, p.\
256]{ON}), thus, if $K$ is a closed Killing vector field, then its
integral curves in $M$ are either circles or points. The converse
of this statement is also true when $(M,g)$ is a compact
Riemannian manifold.
\begin{lem}\label{thm:closedorbitsclosed}
Let $(M,g)$ be a compact Riemannian manifold and let $K$ be a
Killing vector field all of whose integral curves are periodic
(i.e., circles or points). Then $K$ is closed.
\end{lem}
\begin{proof}
Let $G$ be the $1$-parameter group of $\mathfrak x$ and assume by
absurd that $G$ is not closed. Denote by $\overline G$ its closure
in $\mathrm{Iso}(M,g)$, which is a compact abelian Lie group, thus
$\overline G$ is a torus (see for instance \cite[Theorem 10.4]{baker})
of dimension greater than or equal to $2$ (it cannot be $\mathds S^1$,
because there is no $1$-dimensional proper subgroup of $\mathds S^1$).
We claim that the orbits of the actions of $G$ and of
$\overline G$ in $M$ coincide. Namely, if $p\in M$, then clearly
$Gp\subset\overline Gp$; on the other hand, $Gp$ is dense in
$\overline Gp$, because $G$ is dense in $\overline G$. But $Gp$
(and obviously $\overline Gp$) is a closed subset of $M$, because
it is an integral line of $K$, which is a circle or a point, and this proves
that $\overline Gp=Gp$. Since $\overline Gp$ is a circle or a point for all
$p$, it follows that the regular isotropy $H\subset\overline G$ of
the action of the compact group $\overline G$ on $M$ is
non trivial, for otherwise $\overline G$ will be diffeomorphic to a circle or a point.
Since $\overline G$ is abelian, such regular isotropy subgroup is normal in $\overline
G$, which implies that every point of the regular orbits is fixed
by the elements of $H$. This is a contradiction, because the union
of the regular orbits form a dense open subset of $M$ (see for instance
the principal orbit theorem in \cite{Bre}), and no
nontrivial element of $\mathrm{Iso}(M,g)$ acts trivially on a non
empty open subset of $M$.
\end{proof}

\end{proof}

\subsection{Killing vector fields and geodesics}
We will use the following result in \cite[Chapter VI, Proposition~5.7, page 252]{Koba1}, whose simple proof is reproduced here for the reader's convenience.
\begin{lem}\label{thm:intcurvgeo}
Let $(M,g)$ be a semi-Riemannian manifold, and let $K$ be a
Killing vector field on $M$. Let $p_0\in M$ be a critical point
for the function $f(p)=g\big(K_p,K_p\big)$; then, the integral
line of $K$ through $p_0$ is a geodesic.
\end{lem}
\begin{proof}
Since  $f$ is preserved by the flow of the Killing vector $K$,  it
follows that the integral line of $K$ through $p_0$ consists
entirely of critical points of $f$. Thus, it suffices to show
that $(\nabla_KK)_{p_0}=0$. Since $p_0$ is a critical point of $f$
and $K$ is Killing, then for all $v\in T_{p_0}M$ it is:
\[0=v(f)=2g\big(\nabla_vK,K)=-2g(\nabla_KK,v),\]
i.e., $(\nabla_KK)_{p_0}=0$, which concludes the proof.
\end{proof}

\begin{cor}\label{thm:corcurvechiuse}
Let $(M,g)$ be a compact semi-Riemannian manifold with
$\Dim(M)\ge2$ that admits a non trivial Killing vector field $K$
all of whose integral lines are periodic.
Then, there is some non trivial periodic non
self-intersecting geodesic in $M$. If either one of the following two conditions
is satisfied, then there are at least two non trivial
periodic non self-intersecting geodesics in $M$:
\begin{itemize}
\item[(a)] $\min g(K,K)$ and $\max g(K,K)$ are both non zero;
\item[(b)] $K$ is never vanishing.
\end{itemize}
\end{cor}
\begin{proof}
Let $K$ be a Killing vector field as in the assumption. By
Lemma~\ref{thm:intcurvgeo}, the integral curves of $K$ through
critical points of the function $f(p)=g\big(K_p,K_p\big)$ are
geodesics.  Evidently they are periodic; periodic integral curves
of a vector field (considered with minimal period) are non
self-intersecting. If $f$ is identically zero, then every point in
$M$ is critical for $f$, and since $K$ is non trivial then there
are infinitely many non trivial integral curves of $K$ that are
geodesics. If $f$ is not identically zero, then either the minimum
or the maximum of $f$ are non zero, and the corresponding integral
curve of $K$ is a non trivial periodic geodesic.

If either (a) or (b) holds and $f$ is not identically zero, then
the integral curves of $K$ through a minimum and a maximum of $f$
are non trivial distinct periodic geodesics.
\end{proof}

\begin{prop}\label{thm:casoRiemanniano}
Let $(M,g)$ be a compact semi-Riemannian manifold with
$\mathrm{dim}(M)\ge2$ that admits a Killing vector field $K$ that
generates a precompact $1$-parameter subgroup of
$\Iso(M,g)$.  Then $K$ can be
approximated by \emph{closed} Killing vector fields that generate
precompact $1$-parameter subgroups. In particular, there is some
non trivial periodic non self-intersecting geodesic in $M$. If
$K$ satisfies either (a) or (b) of Corollary~\ref{thm:corcurvechiuse}, then there are at least two non
trivial periodic non self-intersecting geodesics in $M$.
\end{prop}
\begin{proof}
Let $\mathfrak x\in\mathfrak{Iso}(M,g)$ be the element
corresponding to $K$, and consider the precompact $1$-parameter
subgroup $G=\big\{\exp(t\mathfrak x):t\in\R\big\}$ of isometries
generated by $K$. Its closure $\overline G$ is a Lie subgroup of
$\mathrm{Iso}(M,g)$, which is abelian and compact, so it is a
torus; denote by $\mathfrak g$ its (abelian) Lie
algebra.  Then, $\mathfrak x$
can be approximated by a sequence $\mathfrak x_n\in\mathfrak g$ of
vectors generating a closed $1$-parameter subgroup of $\overline
G\subset\mathrm{Iso}(M,g)$. The Killing vector fields
$K^n=K^{\mathfrak x_n}$ are closed; using the relation
\eqref{KillingLie} and the compactness of $M$, one sees easily
that
$\lim\limits_{n\to\infty}K^n_p=K_p$ uniformly in $p\in M$.
In particular, if $K$ satisfies either (a) or (b) of Corollary~\ref{thm:corcurvechiuse},
for $n$ large enough also $K^n$ does, and so, the thesis directly follows by
applying Corollary~\ref{thm:corcurvechiuse} to $K^{n}$.
\end{proof}
\subsection{Proof of Theorem}

\begin{prop}\label{thm:elem}
Let $(M,g)$ be a compact Lorentzian manifold and let $K$  be a
Killing vector field on $M$ which is timelike at some
point. Given $H\subset\Iso(M,g)$, assume that for all $\Phi\in H$ and all $q\in M$ it is
$\mathrm d\Phi_{q}(K_{q})=\pm K_{\Phi(q)}$. Then $H$ is precompact.
In particular, the $1$-parameter subgroup of isometries generated
by $K$ is precompact.
\end{prop}
\begin{proof}
Let $p\in M$ be such that $g\big(K_p,K_p\big)<0$. Consider the
compact subsets of $TM$ given by:
\[V=\Big\{\pm K_q:q\in M, \ \text{is such that}\ g\big(K_q,K_q\big)=g\big(K_p,K_p\big)\Big\},\]
and
\[V^\perp=\Big\{v\in K_q^\perp:q\in M\ \text{is such that}\  g\big(K_q,K_q\big)=
g\big(K_p,K_p\big),\ g(v,v)=1\Big\}.\] Consider an orthogonal
basis $b=(v_1,\ldots,v_n)$ of $T_pM$ with $v_1=K_p$ and
$g(v_i,v_j)=\delta_{ij}$ for $i,j\in\{2,\ldots,n\}$. Now recall
that the subset $H\subset\mathrm{Iso}(M,g)$ can be identified with
the $H$-orbit of the basis $b$ by the action of
$\mathrm{Iso}(M,g)$ on the frame bundle $\mathcal F(M)$ (see
\cite[Theorem 1.2, Theorem 1.3]{Koba}). We claim that every vector
of a basis of the $H$-orbit belongs to the compact subset
$V\bigcup V^\perp$, and this implies that the  $H$-orbit of $b$ is
precompact in the frame bundle $\mathcal F(M)$.  The claim follows
easily from the assumption that $\mathrm d\Phi_p\big(K_p\big)=\pm
K_{\Phi(p)}$ for all $\Phi\in H$ and all $p\in M$.

The conclusion applies in particular to the $1$-parameter subgroup $H\subset\mathrm{Iso}(M,g)$
of isometries generated by $K$; for $\Phi\in H$, $\Phi_*K=K$.
\end{proof}

\begin{proof}[Proof of Theorem]
From Proposition~\ref{thm:elem}, the $1$-parameter subgroup of
isometries generated by $K$ is precompact.
Therefore, taking into account that $K$ is timelike
somewhere, the first and second assertions follow directly from
Proposition \ref{thm:casoRiemanniano}. Note that when $K$ is timelike somewhere,
then automatically $\min g(K,K)$ is non zero.

If, in addition, $K$ has at most one periodic integral curve, then
the closed Killing vector fields $K^n$ from Proposition
\ref{thm:casoRiemanniano} that approximate $K$ must be $K^n\ne K$
for all $n$ (recall that $\Dim(M)\ge2$), and we can therefore
assume that the $K^n$'s are pairwise distinct. For each $n$, $K^n$
determines at least two non trivial periodic non self-intersecting
geodesics. The family of all such periodic geodesics for all $n$
cannot be finite. Namely, if it were, then we could find a
subsequence of the $K^n$ with the property that all the elements
of the subsequence have two distinct fixed curves $\gamma_1$ and
$\gamma_2$ as common periodic integral lines. But then, the limit
$K$ would also have $\gamma_1$ and $\gamma_2$ as periodic integral
lines, which gives a contradiction. This concludes the proof.
\end{proof}

Let $(M,g)$ be a stationary Lorentzian manifold and $K$ the
Killing timelike vector field. Consider the auxiliary Riemannian
metric $\gR$ defined using $g$ and $K$ by:
\begin{equation}\label{eq:LorRiem}
\gR(v,w)=g(v,w)-2g\big(v,K_p\big)g\big(w,K_p\big)g\big(K_p,K_p\big)^{-1},
\end{equation}
for all $p\in M$ and $v,w\in T_pM$. The Lorentzian metric $g$ is
given in terms of $\gR$ and $K$ by a similar formula:
\begin{equation}\label{eq:RiemLor}
g(v,w)=\gR(v,w)-2\gR\big(v,K_p\big)\gR\big(w,K_p\big)\gR\big(K_p,K_p\big)^{-1},
\end{equation}
for all $p\in M$ and $v,w\in T_pM$. Since the flow of $K$
preserves the metric $g$ and the field $K$ itself, then using
\eqref{eq:LorRiem} one sees immediately that the flow of $K$
preserves $\gR$, i.e., $K$ is a Killing vector field also for the
metric $\gR$. Using this fact, we get the following immediate
corollary of Lemma~\ref{thm:closedorbitsclosed}:
\begin{cor}
Let $(M,g)$ be a compact stationary Lorentzian manifold, and let
$K$ be a timelike Killing vector field on $M$. Then, $K$ is closed
if and only if all its integral lines are periodic.\qed
\end{cor}
\begin{rem}
Even when the two geodesics $\gamma_1$ and $\gamma_2$ determined in
our main theorem (as integral lines of a closed vector field) are
timelike, it does not necessarily follow that they belong to the
same free (timelike) homotopy class, see
Remark~\ref{rem:kleinbottle}, but rather that they have some
\emph{iterate} that belong to the same free (timelike) homotopy
class.
\end{rem}

\end{section}

\begin{section}{On the topological structure\\ of a compact
stationary Lorentzian manifold}\label{sec:topology}
\subsection{Fibration associated to a closed Killing vector field}\label{sub:fibration}
A compact manifold $M$ will be called a \emph{generalized Seifert fibered
space} if it admits a smooth action of the circle $\mathds S^1$ without
fixed points or, equivalently, with finite isotropy. The orbits of a fixed point free
action of $\mathds S^1$, that are diffeomorphic to $\mathds S^1$, are called the
\emph{fibers} of the fibered space. Low dimensional generalized
Seifert fibered spaces are classified, see \cite{Fin1, Fin2, Kol, Sco, Sei}.
By standard results on group actions, probably going back to Seifert,
the orbit space of a smooth action of a compact Lie group on a compact manifold
having finite isotropy has the structure of a compact orbifold (see the Appendix of E.\ Salem in \cite{Molino}
for details on orbifolds; the book contains also a more general result on the orbifold structure of
orbit spaces in the context of Riemannian foliations).

Using the results of Section~\ref{sec:proofoftheorem}, it is easy to prove the following:
\begin{prop}\label{thm:propfibration}
A compact manifold $M$ admits a Lorentzian metric tensor with a
never vanishing Killing vector field that is timelike somewhere
if and only if it is diffeomorphic to a
generalized Seifert fibered space.  In this case, the metric can be chosen to
have a timelike Killing vector field.
\end{prop}
\begin{proof}
As shown in the proof of the Theorem, a Lorentzian manifold as in
the hypotheses above has a never vanishing closed Killing vector
field $K$. The one-parameter group of isometries generated by such
a Killing field gives a smooth action of $\mathds S^1$ without
fixed points; $K$ is tangent to the fibers of this action.
Conversely, given a smooth action of $\mathds S^1$ on $M$ without
fixed points, by a standard averaging argument one can find a
Riemannian metric tensor $\gR$ which makes such action isometric,
i.e., the infinitesimal generator $K$ of this action is
$\gR$-Killing (see for instance \cite{knapp}). Consider the
Lorentzian metric tensor $g$ defined as in \eqref{eq:RiemLor};
then $K$ is timelike and $g$-Killing.
%
\end{proof}

\begin{rem}\label{rem:kleinbottle}
Given a \emph{free} action of $\mathds S^1$ on a compact manifold
$M$, then the orbit space $M/\mathds S^1$ is a smooth manifold
(see for instance \cite[Theorem 23.4]{cannas} or \cite[Theorem
4.3]{Koba1}). We observe however that in general the quotient
space $M_0=M/\mathds S^1$ is not a manifold. As an example,
consider $M$ to be the Lorentzian Klein bottle obtained as the
quotient of $\R^2$ endowed with the Minowski metric $\mathrm
dx^2-\mathrm dt^2$ by the action of the group generated by the
isometries $(x,t)\mapsto(x+1,t)$ and $(x,t)\mapsto(1-x,t+1)$. The
vector field $K=\frac\partial{\partial t}$ on $M$ is timelike and
Killing; all its integral lines are periodic. It is easily seen
that in this case the $\mathds S^1$-action induced by the flow of
$K$ has exactly two exceptional orbits, and that the orbit space
$M/\mathds S^1$ is homeomorphic to the closed interval
$[0,\frac12]$. Note that the periodic integral lines of $K$
corresponding to the two exceptional orbits do \emph{not} belong
to the same free homotopy class of the other integral lines of
$K$, but rather their two-fold iteration is in the free homotopy
class of the other integral lines of $K$.
\end{rem}

As a corollary of Proposition~\ref{thm:propfibration}, we get a somewhat better estimate on the number of
periodic geodesics given in terms of the
\emph{Ljusternik--Schnirel\-man category}. Recall that the
Ljusternik--Schnirel\-man category (shortly, LS category)
$\mathrm{cat}(\mathcal X)$ of a
topological space $\mathcal X$ is the cardinality (possibly
infinite) of a minimal family of closed contractible subsets of
$\mathcal X$ whose union covers $\mathcal X$. If $\mathcal X$ is $G$-space, i.e., a
topological space on which a compact  group $G$ is acting
continuously, then one can define the equivariant notion of
Ljusternik--Schnirel\-man $G$-category $\mathrm{cat}_{G}(\mathcal
X)$ (see for instance \cite{Marz}).  A homotopy $H:U\times[0,1]\to\mathcal X$
of an open $G$-invariant set $U\subset\mathcal X$ is called \emph{$G$-equivariant}
if $gH(x,t)=H(gx,t)$ for any $g\in G$, $x\in U$ and $t\in[0,1]$.
The set $U$ is \emph{$G$-categorical} if there is a $G$-homotopy
$H$ with $H(\cdot,0)$ the identity, and $H(\cdot,1)$ maps $U$ to a single orbit.
The \emph{equivariant category} $\mathrm{cat}_{G}(\mathcal X)$ is the cardinality
of a minimal family of $G$-categorical
open sets whose union covers $\mathcal X$.

If $G$ is a compact Lie group, $\mathcal X$ is a smooth $G$-manifold,
and  $h:\mathcal X\to\R$ is a smooth function which is $G$-invariant, then $h$ has at
least $\mathrm{cat}_{G}(\mathcal X)$ distinct critical $G$-orbits (see
\cite[Th.~3.2]{Marz}).
\begin{cor}\label{thm:LS}
Let $(M,g)$ be a compact Lorentzian manifold with a never
vanishing closed Killing vector field on $M$.
Consider the $\mathds S^1$-action on $M$
determined by $K$. Then, there are at least $\mathrm{cat}_{\mathds
S^1}(M)$ distinct periodic non self-intersecting geodesics in $M$.
\end{cor}
\begin{proof}
The function $f:M\to\R$ defined by $f(p)=g\big(K_p,K_p\big)$ is
constant on the orbits of $G=\mathds S^1$, thus it has at least
$\mathrm{cat}_{G}(M)$ critical orbits. Hence, the proof follows by
observing that distinct critical $G$-orbits of $f$ in $M$
correspond to distinct non self-intersecting periodic geodesics.
\end{proof}

In Corollary~\ref{thm:LS}, note that $\mathrm{cat}_{\mathds S^1}(M)\ge2$.
Namely, if it were $\mathrm{cat}_{\mathds S^1}(M)=1$, then $M$ would be
(equivariantly) homotopic to an orbit of $\mathds S^1$, which is diffeomorphic to $\mathds S^1$.
But, no compact manifold of dimension greater than or equal to $2$ is
homotopic to $\mathds S^1$. Observe also that in general the equivariant
LS category $\mathrm{cat}_{\mathds S^1}(M)$ is greater than or equal to the
LS category $\mathrm{cat}(M/\mathds S^1)$ of the quotient space $M/\mathds S^1$.
The Klein bottle in Remark~\ref{rem:kleinbottle} provides an example where such
inequality is strict: here the quotient space $M/\mathds S^1$ is contractible,
and thus $\mathrm{cat}(M/\mathds S^1)=1$, while it is easily computed $\mathrm{cat}_{\mathds S^1}(M)=2$.

\subsection{Some examples}\label{sub:someexamples}
It is easy to produce examples of compact stationary Lorentzian
manifolds with a timelike Killing vector field having all of its
integral curves periodic or having no periodic integral curve at
all. In next examples we show two different constructions to
produce compact stationary manifolds with a timelike Killing
vector field having integral curves of mixed type.
\begin{example}\label{exa:dueorbitechiuse}
Consider the following smooth isometric action of $\mathds T^2=\mathds
S^1\times\mathds S^1$ on the round $3$-sphere $\mathds S^3$. Set
$\mathds S^3=\big\{(z,w)\in\C^2:\vert z\vert^2+\vert
w\vert^2=1\big\}$, and for $(\lambda_1,\lambda_2)\in\mathds
S^1\times\mathds S^1$, $(z,w)\in\mathds S^3$, define
$(\lambda_1,\lambda_2)\cdot(z,w)=(\lambda_1z,\lambda_2w)\in\mathds
S^3$. The circles $c_1=\mathds
S^1\times\{0\}\subset\mathds S^3$ and
 $c_2=\{0\}\times\mathds S^1\subset\mathds S^3$ are orbits of this action, and they have non trivial
 isotropy. Every other orbit is regular, and it is diffeomorphic to $\mathds T^2$.
 Now consider a $1$-parameter subgroup $G\subset\mathds T^2$ which is dense in $\mathds T^2$
 and the restriction of the action of $\mathds T^2$ to $G$. The singular orbits $c_1$ and $c_2$ are also
 orbits of $G$, since the projections $\pi_1,\pi_2:G\to\mathds S^1$ given by the restrictions
 of the projections $\pi_1,\pi_2:\mathds S^1\times\mathds S^1\to\mathds S^1$ onto the first and the second
 factor respectively, are surjective. All the other orbits of $G$ are clearly not closed.
Consider the Killing vector field $K$ for the Riemannian metric
$\gR$ corresponding to the isometric action of $G$, and define a
Lorentzian metric $g$ on $\mathds S^3$ using formula
\eqref{eq:RiemLor}. Then, $K$ is a timelike Killing vector field
for $g$, and we obtain an example of a compact stationary
Lorentzian manifold admitting a timelike Killing vector field with
exactly two periodic integral lines.
\end{example}
Next, we show a construction of an isometric action
of $\R$ on a compact manifold with only one periodic orbit.
\begin{example}\label{exa:unaorbitachiusa}
Let $(M_0,h)$ be a compact Riemannian manifold (without boundary)
that admits an isometry $\psi:M_0\to M_0$ which has exactly one
fixed point $p_0$ and having no other periodic point, i.e., all
the powers $\psi^N$, $N\in\N\setminus\{0\}$, only have $p_0$ as
fixed point.

An example of this situation can be obtained as follows. Consider
the rotation $R_\theta:\mathds S^2\to\mathds S^2$ around the
north-south axis by an angle $\theta\in\left[0,2\pi\right[$ which
is not a rational multiple of $\pi$. Here the two-sphere $\mathds
S^2$ is endowed with the round metric. Then, $R_\theta$ induces
a map $\overline R_\theta:\R\mathrm P^2\to\R\mathrm P^2$ on the
projective plane $\R\mathrm P^2$ which is an isometry with only one fixed point
and no other periodic point.

Let $M$ be the manifold obtained as a quotient of the product
$M_0\times\R$, by identifying $(p,n)$ with
$\big(\psi(p),n+1\big)$, $n\in\mathds Z$. The Lorentzian metric
$h\oplus(-\mathrm dt^2)$ on $M_0\times\R$ induces a Lorentzian
metric $g$ on $M$, for which the maps
$T_s$ induced by the translations $M_0\times\R\ni(p,t)\to(p,s+t)\in M_0\times\R$
are isometries.  The $1$-parameter group
of isometries of $(M,g)$ given by $\R\ni s\mapsto
T_s\in\mathrm{Iso}(M,g)$ has a timelike Killing vector field
$K=\frac\partial{\partial t}$ as infinitesimal generator. It is
easy to see that $K$ has exactly one closed integral curve, which
is the one passing through $p_0$.
\end{example}
\subsection{On the isometry group of a compact stationary Lorentzian manifold}\label{sub:isometrygroup}
As to the isometry group of a compact stationary Lorentzian manifold, it is
known that it may fail to be compact (see for instance
\cite[Remark 4.3]{San0} or \cite{Da}).
In fact, there exists a complete classification of Lie groups that appear as
connected components of the identity of the isometry group of compact Lorentzian manifolds,
which is due independently to Adams/Stuck \cite{AdamsStuck} and to Zeghib \cite{Ghani1}.

We have the following partial results concerning the compactness of the isometry group.
\begin{prop}
Let $(M,g)$ be a compact Lorentzian manifold with a never
vanishing closed Killing vector field $K$ that is timelike
somewhere.
 Let $M_0=M/\mathds S^1$ be the orbit space of the
corresponding $\mathds S^1$-action and let $\pi:M\to M_0$ be the
canonical projection. Then, the group $\mathrm{Iso}(M,g;\pi)$
consisting of all
isometries of $(M,g)$ that preserve the fibration is compact.
\end{prop}
\begin{proof}
 A diffeomorphism $\Phi$ of $M$ preserves the
fibration if and only if the push-forward $\Phi_*(K)$ is a
pointwise multiple of $K$, that is, $\Phi_{*}(K)=\lambda K$ for
some function $\lambda$ on $M$. Moreover, as $\Phi$ is an isometry,
$\Phi_{*}K$ must be Killing and it can be easily proved that
$\lambda$ has to be constant (if $K$ is Killing, $\lambda K$ is
Killing iff $\lambda$ is constant). Moreover,  we have that
\begin{multline*}\min_{p\in M}g(K_p,K_p)=\min_{p\in M}g\big(K_{\Phi^{-1}(p)},K_{\Phi^{-1}(p)}\big)=
\min_{p\in M}g\big((\Phi_*K)_{p},(\Phi_*K)_{p}\big)\\=\lambda^2\min_{p\in M}g(K_p,K_p),\end{multline*}
and hence $\lambda^2=1$.
Applying Proposition~\ref{thm:elem} to
$H=\mathrm{Iso}(M,g;\pi)$ and $K$, we obtain that
$\mathrm{Iso}(M,g;\pi)$ is precompact. But it is clearly closed.
Whence, $\mathrm{Iso}(M,g;\pi)$ is compact.
%
\end{proof}
It is proven in \cite[Lemma 4.4]{San0} that a compact Lorentzian manifold
admitting a Killing vector field which is timelike at one point
and whose isometry group is one-dimensional, then the isometry
group must be compact. We have the following similar result:
\begin{cor}
Let $(M,g)$ be a compact Lorentzian manifold with a Killing vector
field that is timelike somewhere and assume that
$\mathrm{Iso}(M,g)$ is abelian. Then, $\mathrm{Iso}(M,g)$  is compact.
\end{cor}
\begin{proof} If $\mathrm{Iso}(M,g)$ is abelian then $\mathrm{Ad}_{\Phi}=\mathrm{Id}$ for all $\Phi$. Therefore,
the thesis directly follows from Lemma~\ref{thm:Adpushforward} and Proposition~\ref{thm:elem} applied to $H=\mathrm{Iso}(M,g)$.
\end{proof}
\end{section}

\begin{section}{Final results and remarks}\label{sec:final}

The technique used in this paper for proving the existence of
periodic geodesics is not limited to the hypotheses of our main
theorem. It can essentially be applied to any compact
semi-Riemannian manifold $(M,g)$ under the more general hypotheses
of Proposition \ref{thm:casoRiemanniano}. To illustrate this, we
give below some simple examples:
\begin{example}
Any compact
two-dimensional Lorentzian manifold $(M,g)$ admits some timelike
periodic geodesic \cite{Gal2}. If, in addition, we assume that
$M\cong T^2$ is diffeomorphic to a torus and that it admits a
(non-trivial) Killing vector field $K$, then $K$ does not vanish
at any point \cite[Th. 4.2]{San1}. Moreover, we have two
possibilities: either $(T^{2},g)$ is flat, and so, it contains
infinitely many periodic timelike geodesics; or $(T^{2},g)$ is non
flat, which from \cite[Th. 4.2]{San1} implies that Iso$(T^{2},g)$
is compact. Whence, applying Proposition
\ref{thm:casoRiemanniano}:
\begin{prop} Any Lorentzian torus $(T^{2},g)$ with a Killing vector field $K\not\equiv 0$
contains at least two \emph{geometrically distinct} non trivial
periodic non self-intersecting geodesics.\qed
\end{prop}
\end{example}
\begin{example}
From \cite{Da}, any simply connected compact
real-analytic Lorentzian manifold $(M,g)$ has compact isometry
group. Therefore:
\begin{prop} Any simply connected compact real-analytic Lorentzian manifold $(M,g)$
with $\mathrm{dim}\big(\mathrm{Iso}(M,g)\big)>0$ contains some non trivial
periodic geodesic.
\end{prop}
\end{example}
\begin{example}
 Consider a compact
semi-Riemannian manifold $(M,g)$ of index $m\ge1$ that admits $m$
Killing vector fields $K^1,\ldots
K^{m}$ which generate a negative definite subspace for $g$ of
dimension $m$ at some point $p$. If $m=1$, then we are in the Lorentzian case studied in this paper.
The geodesic connectedness in the particular case of
\emph{generalized stationary semi-Riemannian manifolds} (i.e., semi-Riemannian
manifolds that admit a timelike distribution $D\subset TM$ of rank $m$ which is
generated by $m$ pointwise linearly independent commuting timelike
Killing vector fields $K^{1},\ldots, K^{m}$) is studied in \cite{GiaPicSam}.
Denote by $\mathcal A:M\to\mathrm{Mat}_m(\R)$ the map taking values in the space of $m\times m$ real
symmetric matrices
defined by $\mathcal A(q)=\big(g(K^i_q,K^j_q)\big)_{ij}$, and
consider the compact subsets of $TM$ given by:
\[V=\Big\{ K^i_q: i=1,\ldots,m,\ q\in M\ \text{is such that}\ \mathcal A(q)=\mathcal A(p)\Big\}\]
and
\[V^\perp=\Big\{v\in \bigcap_{i=1}^m(K^i_q)^\perp:q\in M\ \text{is as above and}\  g(v,v)=1\Big\}.\]
Consider a basis $b=(v_1,\ldots,v_n)$ of $T_pM$ with $v_i=K^{i}_p$
for all $1\leq i\leq m$ and $g(v_i,v_j)=\delta_{ij}$ for
$i\in\{m+1,\ldots,n\}$ and all $j$. The Killing vector fields
$K^1\ldots, K^m$ generate a precompact subgroup $H$ of
$\mathrm{Iso}(M,g)$ (as in Proposition~\ref{thm:elem}). Under the
assumptions that the $K^i$ commute, i.e., $[K^i,K^j]=0$ for all
$i,j=1,\ldots,m$, then $H$ is abelian, and thus the closure
$\overline H$ is a torus. Consider a non zero vector $\mathfrak x$
in the Lie algebra of $\overline H$ such that the corresponding
Killing vector field $K^\mathfrak x$ is closed. Then, as in the
proof of our main result, there is some integral curve of
$K^\mathfrak x$ that is a geodesic, say $\gamma:\R\to M$. If the
$K^i$'s generate a distribution of rank $k>1$, we claim that there
are infinitely many periodic integral curves of $K^\mathfrak x$
that are geodesics. Namely, if $k>1$, there is at least one
$l\in\{1,\ldots,m\}$ such that $K^l_{\gamma(0)}$ and $K^\mathfrak
x_{\gamma(0)}$ are linearly independent. Denote by
$(\varphi_t)_{t\in\R}$ the flow of $K^l$; since $K^l$ and
$K^{\mathfrak x}$ commute then, for all $t\in\R$,
$\varphi_t\circ\gamma$ is an integral line of $K^\mathfrak x$. On
the other hand, since $K^l$ is Killing, $\varphi_t\circ\gamma$ is
a periodic geodesic in $(M,g)$. Since $K^l_{\gamma(0)}$ and
$K^\mathfrak x_{\gamma(0)}$ are linearly independent, then for
$t\in\R$ sufficiently small the curves $\varphi_t\circ\gamma$ and
$\gamma$ are distinct. We have proven the following:
\begin{prop} Compact semi-Riemannian manifolds as above of index $k>1$ have infinitely
many distinct periodic geodesics.
\end{prop}
\end{example}
\end{section}

\end{document}